\documentclass[a4paper,12pt]{amsart}
\usepackage{graphics}
\usepackage{amssymb}
\usepackage{amsmath}
\usepackage{latexsym}
\usepackage{amscd}
\newtheorem{Th}{Theorem}
\newtheorem{Prop}[Th]{Proposition}

\newtheorem{Rem}[Th]{\sc Remark}

\newtheorem{Lemma}[Th]{Lemma}
\newcommand{\be}{\begin{eqnarray*}}
\newcommand{\ee}{\end{eqnarray*}}
\newcommand{\pkEF}{{\mathcal P}(^k\!E,F)}

\newcommand{\LEF}{{\mathcal L}(E,F)}

\newcommand{\HEF}{{\mathcal H}(E,F)}
\newcommand{\HUEF}{{\mathcal H}_{\mathcal U}(E,F)}

\newcommand{\pUkEF}{{\mathcal P}_{\mathcal U}(^k\!E,F)}
\newcommand{\CU}{{\mathcal C}_{\mathcal U}}

\newcommand{\Proof}{\noindent {\sc Proof. }}

\newcommand{\finesp}{\hspace*{\fill} $\Box$\vspace{.5\baselineskip}}

\newcommand{\N}{\ensuremath{\mathbb{N}}}

\newcommand{\ra}{\rightarrow}
\newcommand{\lra}{\longrightarrow}
\newcommand{\Ra}{\Rightarrow}       
\newcommand{\LRa}{\Longrightarrow}  
\newcommand{\espv}{\vspace{.5\baselineskip}}
\newcommand{\espvv}{\vspace{\baselineskip}}

\newcommand{\eps}{\epsilon}
\newcommand{\kdots}{\stackrel{(k)}{\ldots}}

\begin{document}
\title{Surjective factorization of holomorphic mappings}

\author[M. Gonz\'alez]{Manuel Gonz\'alez}
\address{Departamento de Matem\'aticas \\
      Facultad de Ciencias\\
      Universidad de Cantabria \\ 39071 Santander (Spain)}
\email{gonzalem@ccaix3.unican.es}
\thanks{The first named author was supported
in part by DGICYT Grant PB 97--0349 (Spain)}

\author[J. M. Guti\'errez]{Joaqu\'\i n M. Guti\'errez}
\address{Departamento de Matem\'atica Aplicada\\
      ETS de Ingenieros Industriales \\
      Universidad Polit\'ecnica de Madrid\\
      C. Jos\'e Guti\'errez Abascal 2 \\
      28006 Madrid (Spain)}
\email{jgutierrez@etsii.upm.es}
\thanks{The second named author was supported
in part by DGICYT Grant PB 96--0607 (Spain)}
\thanks{\hspace*{\fill}\scriptsize file sfhm.tex}



\begin{abstract}
We characterize the holomorphic mappings $f$ between complex
Ba\-nach spaces that may be written in the form $f=T\circ g$,
where $g$ is another holomorphic mapping and $T$ belongs to a
closed surjective operator ideal.
\end{abstract}

\maketitle

\section{Introduction and preliminary results}

In recent years many authors \cite{AG,AS,Ge,GG,GGiny,Li,Ro,RyWCo}
have studied conditions on a holomorphic mapping $f$ between
complex Banach spaces so that it may be written in the form either
$f=g\circ T$ or $f=T\circ g$, where $g$ is another holomorphic
mapping and $T$ a (linear bounded) operator belonging to certain
classes of operators.

A rather thorough study of the factorization of the form $f=g\circ
T$, where $T$ is in a closed injective operator ideal, was carried
out by the authors in \cite{GGiny}. In the present paper we
analyze the case $f=T\circ g$.

If $f=T\circ g$, with $T$ in the ideal of compact operators, and
$g$ is holomorphic on a Banach space $E$ then, since $g$ is
locally bounded, $f$ will be ``locally compact" in the sense that
every $x\in E$ has a neighbourhood $V_x$ such that $f(V_x)$ is
relatively compact. It is proved in
\cite{AS} that the converse also holds: every locally compact
holomorphic mapping $f$ can be written in the form $f=T\circ g$,
with $T$ a compact operator. Similar results were given in
\cite{RyWCo} for the ideal of weakly compact operators,
in \cite{Li} for the Rosenthal operators, and in \cite{Ro} for the
Asplund operators. We extend this type of factorization to every
closed surjective operator ideal.\espv

Throughout, $E$, $F$ and $G$ will denote complex Banach spaces,
and \N\ will be the set of natural numbers. We use $B_E$ for the
closed unit ball of $E$, and $B(x,r)$ for the open ball of radius
$r$ centered at $x$. If $A\subset E$, then $\bar{\Gamma}(A)$
denotes the absolutely convex, closed hull of $A$, and if $f$ is a
mapping on $E$, then
$$
\|f\|_A:=\sup\{ |f(x)|:x\in A\}\, .
$$

We denote by $\LEF$ the space of all operators from $E$ into $F$,
endowed with the usual operator norm. A mapping $P:E\ra F$ is a
$k$-homogeneous (continuous) {\it polynomial\/} if there is a
$k$-linear continuous mapping $A:E\times\kdots\times E\ra F$ such
that $P(x)=A(x,\ldots,x)$ for all $x\in E$. The space of all such
polynomials is denoted by $\pkEF$. A mapping $f:E\ra F$ is {\it
holomorphic\/} if, for each $x\in E$, there are $r>0$ and a
sequence $(P_k)$ with $P_k\in\pkEF$ such that
$$
f(y)=\sum_{k=0}^\infty P_k(y-x)
$$
uniformly for $\|y-x\|<r$. We use the notation
$$
P_k=\frac{1}{k!}\,d^k\! f(x),
$$
while $\HEF$ stands for the space of all holomorphic mappings from
$E$ into $F$.

We say that a subset $A\subset E$ is {\it circled\/} if for every
$x\in A$ and complex $\lambda$ with $|\lambda |=1$, we have
$\lambda x\in A$.

For a general introduction to polynomials and holomorphic
mappings, the reader is referred to \cite{Di,Mu,Na}. The
definition and general properties of operator ideals may be seen
in
\cite{Pi}.

An operator ideal ${\mathcal U}$ is said to be {\it injective\/}
\cite[4.6.9]{Pi} if, given an operator $T\in\LEF$ and an injective
isomorphism $i:F\ra G$, we have that $T\in{\mathcal U}$ whenever
$iT\in{\mathcal U}$. The ideal ${\mathcal U}$ is {\it
surjective\/}
\cite[4.7.9]{Pi} if, given $T\in\LEF$ and a surjective operator
$q:G\ra E$, we have that $T\in {\mathcal U}$ whenever
$Tq\in{\mathcal U}$. We say that ${\mathcal U}$ is {\it closed\/}
\cite[4.2.4]{Pi} if for all $E$ and $F$, the space ${\mathcal
U}(E,F):=\{ T\in\LEF :T\in{\mathcal U}\}$ is closed in $\LEF$.

Given an operator $T\in\LEF$, a procedure is described in
\cite{DFJP} to construct a Banach space $Y$ and operators
$k\in{\mathcal L}(E,Y)$ and $j\in{\mathcal L}(Y,F)$ so that
$T=jk$. We shall refer to this construction as the DFJP {\it
factorization.} It is shown in \cite[Propositions~1.6 and~1.7]{H}
(see also
\cite[Proposition~2.2]{Go} for simple statement and proof) that
given an operator $T\in\LEF$ and a closed operator ideal
${\mathcal U}$,

(a) if ${\mathcal U}$ is injective and $T\in {\mathcal U}$, then
$k\in{\mathcal U}$;

(b) if ${\mathcal U}$ is surjective and $T\in {\mathcal U}$, then
$j\in{\mathcal U}$.\\ We say that ${\mathcal U}$ is {\it
factorizable\/} if, for every $T\in{\mathcal U}(E,F)$, there are a
Banach space $Y$ and operators $k\in{\mathcal L}(E,Y)$ and
$j\in{\mathcal L}(Y,F)$ so that $T=jk$ and the identity $I_Y$ of
the space $Y$ belongs to ${\mathcal U}$.

We now give a list of closed operator ideals which are injective,
surjective or factorizable. We recall the definition of the most
commonly used, and give a reference for the others.

An operator $T\in\LEF$ is {\it (weakly) compact\/} if $T(B_E)$ is
a relatively (weakly) compact subset of $F$; $T$ is {\it (weakly)
completely continuous\/} if it takes weak Cauchy sequences in $E$
into (weakly) convergent sequences in $F$; $T$ is {\it
Rosenthal\/} if every sequence in $T(B_E)$ has a weak Cauchy
subsequence; $T$ is {\it unconditionally converging\/} if it takes
weakly unconditionally Cauchy series in $E$ into unconditionally
convergent series in $F$.\espv

\begin{tabular}{|l|c|c|c|} \hline\hline\espv
{\it Closed operator ideals} & {\it Injective} & {\it Surjective} &
                                {\it Factorizable}\espv\\ \hline\espv
compact operators & {\bf Yes} & {\bf Yes} & No
\\ weakly compact  & {\bf Yes} & {\bf Yes}  & {\bf
Yes} \\ Rosenthal & {\bf Yes} & {\bf Yes} & {\bf Yes}
\\ completely continuous & {\bf Yes} & No & No
\\ weakly completely continuous & {\bf Yes} & No &
No \\ unconditionally converging & {\bf Yes} & No
& No \\ Banach-Saks \cite[\S3]{JaWCO} & {\bf Yes} & {\bf Yes} &
{\bf Yes} \\ weakly Banach-Saks \cite[\S3]{JaWCO} & {\bf Yes} & No
& No
\\ strictly singular \cite[1.9]{Pi} & {\bf Yes} & No
& No \\ separable range   & {\bf Yes} & {\bf Yes} & {\bf Yes} \\
strictly cosingular \cite[1.10]{Pi} & No & {\bf Yes} & No \\
limited \cite{BD}   & No & {\bf Yes} & No \\ Grothendieck
\cite{DLS}  & No & {\bf Yes} & No \\ decomposing (Asplund)
\cite[24.4]{Pi}  & {\bf Yes} & {\bf Yes} & {\bf Yes} \\
Radon-Nikod\'ym \cite[24.2]{Pi} & {\bf Yes} & No & No \\
absolutely continuous \cite[\S3]{JM} & {\bf Yes} & No & No\espv\\
\hline\hline
\end{tabular}\espvv

The results on this list may be found in \cite{Pi} and the other
references given, for the injective and surjective case. The
factorizable case may be seen in \cite{H}.

If ${\mathcal U}$ is an operator ideal, the {\it dual ideal\/}
${\mathcal U}^d$ is the ideal of all operators $T$ such that the
adjoint $T^*$ belongs to ${\mathcal U}$. Easily, we have:
\be
{\mathcal U} \mbox{ is closed injective } & \LRa & {\mathcal U}^d
\mbox{ is closed surjective }\\
{\mathcal U} \mbox{ is closed surjective } & \LRa & {\mathcal U}^d
\mbox{ is closed injective }
\ee
The list above might therefore be completed with some more dual
ideals.

Moreover, to each $T\in\LEF$ we can associate an operator
$T^q:E^{**}/E\to F^{**}/F$ given by
$T^q(x^{**}+E)=T^{**}(x^{**})+F$. Let ${\mathcal U}^q:=\{
T\in\LEF: T^q\in{\mathcal U}\}$. Then, if ${\mathcal U}$ is
injective (resp.\ surjective, closed), so is ${\mathcal U}^q$
\cite[Theorem~1.6]{Go}.

\begin{Rem}{\rm
There is another notion of factorizable operator ideal which may
be used. We say that ${\mathcal U}$ is DFJP {\it factorizable\/}
\cite[Definition~2.3]{Go} if, for every $T\in{\mathcal U}$, the
identity of the intermediate space in the DFJP factorization of
$T$ belongs to ${\mathcal U}$. Clearly, every DFJP factorizable
operator ideal is factorizable. The following example shows that
the converse is not true. Let ${\mathcal A}$ be the ideal of all
the operators that factor through a subspace of $c_0$. Clearly,
${\mathcal A}$ is factorizable. Consider the operator
$T:\ell_2\ra\ell_2$ given by $T((x_n)):=(x_n/n)$. We have
$T\in{\mathcal A}$. The intermediate space in the DFJP
factorization is an infinite dimensional reflexive space. Clearly,
the identity map on it does not belong to ${\mathcal A}$.

All the factorizable ideals on the table above are DFJP
factorizable \cite{Go}. Note also that, if ${\mathcal U}$ is DFJP
factorizable, then so are ${\mathcal U}^d$ and ${\mathcal U}^q$
\cite{Go}.}
\end{Rem}

\section{Surjective factorization}

In this Section, we study the factorizations in the form $T\circ
g$, with $T\in{\mathcal U}$, where $\mathcal U$ is a closed
surjective operator ideal.

\begin{Lemma}
\label{surjid}
{\rm\cite[Proposition~2.9]{JaWCO}} Given a closed surjective
operator ideal $\mathcal U$, let $S\in\LEF$ and suppose that for
every $\eps>0$ there are a Banach space $D_\eps$ and an operator
$T_\eps\in{\mathcal U}(D_\eps,F)$ such that
$$
S(B_E)\subseteq T_\eps (B_{D_\eps}) +\eps B_F\, .
$$
Then, $S\in{\mathcal U}$.
\end{Lemma}

We denote by $\CU(E)$ the collection  of all $A\subset E$ so that
$A\subseteq T(B_Z)$ for some Banach space $Z$ and some operator
$T\in{\mathcal U}(Z,E)$ (see \cite{S}).

The following probably well-known properties of $\CU$ will be
needed:

\begin{Prop}
\label{propCU}
Let $\mathcal U$ be a closed surjective operator ideal. Then:

{\rm (a)} If $A\in \CU(E)$ and $B\subset A$, then $B\in\CU(E)$;

{\rm (b)} if $A_1,\ldots,A_n\in\CU(E)$, then
$\cup_{i=1}^nA_i\in\CU(E)$ and $\sum_{i=1}^nA_i\in\CU(E)$;

{\rm (c)} if $A\subset E$ is bounded and, for every $\eps>0$,
there is a set $A_\eps\in\CU(E)$ such that $A\subseteq A_\eps
+\eps B_E$, then $A\in\CU(E)$.

{\rm (d)} if  $A\in \CU(E)$, then $\bar{\Gamma}(A)\in\CU(E)$;
\end{Prop}

\Proof
(a) is trivial and (b) is easy. Both are true without any
assumption on the operator ideal ${\mathcal U}$.

(c) For $A\subset E$ bounded, consider the operator
$$
T:\ell_1(A)\lra E\quad\mbox{ given by }\quad T\left(
(\lambda_x)_{x\in A}\right)=\sum_{x\in A}\lambda_xx.
$$
Given $\eps>0$, there is $A_\eps\in\CU(E)$ such that $A\subseteq
A_\eps +\eps B_E$. Therefore,
$$
A\subseteq T\left(
B_{\ell_1(A)}\right)\subseteq\bar{\Gamma}(A)\subseteq
\Gamma (A)+\eps B_E\subseteq \Gamma (A_\eps)+2\eps B_E.
$$
Clearly, $\Gamma (A_\eps)\in\CU(E)$. Hence, $T\in{\mathcal U}$ (by
Lemma~\ref{surjid}), and $A\in\CU(E)$.

(d) If $A\in\CU(E)$, there is a space $Z$ and $T\in{\mathcal
U}(Z,E)$ such that $A\subseteq T(B_Z)$. Therefore, for all
$\eps>0$,
$$
\bar{\Gamma}(A)\subseteq \overline{T(B_Z)}\subseteq T(B_Z)+\eps
B_E.
$$
Now, it is enough to apply part (c).\finesp

We shall denote by $\HUEF$ the space of all $f\in\HEF$ such that
each $x\in E$ has a neighbourhood $V_x$ with $f(V_x)\in \CU(F)$.
Easily, a polynomial $P\in\pkEF$ belongs to $\HUEF$ if and only if
$P(B_E)\in\CU(F)$. The set of all such polynomials will be denoted
by ${\mathcal P}_{\mathcal U}(^k\!E,F)$.

The following result is an easy consequence of the Hahn-Banach
theorem and the Cauchy inequality

\begin{Lemma}
\label{Cauchy}
{\rm \cite[Lemma~3.1]{RyWCo}}
Given $f\in\HEF$, a circled subset $U\subset E$, and $x\in E$, we have
$$
\frac{1}{k!}\,d^k\! f(x)(U)\subseteq \bar{\Gamma}(f(x+U))
$$
for every $k\in\N$.
\end{Lemma}

\begin{Prop}
Let $\mathcal U$ be a closed surjective operator ideal, and
$f\in\HEF$. The following assertions are equivalent:

{\rm (a)} $f\in\HUEF$;

{\rm (b)} there is a zero neighbourhood $V\subset E$ such that
$f(V)\in \CU(F)$;

{\rm (c)} for every $k\in\N$ and every $x\in E$, we have that
$d^k\! f(x)\in\pUkEF$;

{\rm (d)} for every $k\in\N$, we have that $d^k\! f(0)\in\pUkEF$.
\end{Prop}

\Proof
(a) $\Ra$ (c) and (b) $\Ra$ (d) follow from Lemma~\ref{Cauchy}.

(d) $\Ra$ (a) Let $x\in E$. There is $\eps>0$ such that
$$
f(y)=\sum_{k=0}^\infty \frac{1}{k!} \, d^k\!f(0)(y)
$$
uniformly for $y\in B(x,\eps)$ \cite[\S7, Proposition~1]{Na}.
By Proposition~\ref{propCU}(b), for each $m\in \N$, we have
$$
\left\{ \sum_{k=0}^m \frac{1}{k!}\,d^k\!f(0)(y):y\in B(x,\eps)\right\}
\in \CU(F)\, .
$$
Using the uniform convergence on $B(x,\eps)$, and
Proposition~\ref{propCU}(c), we conclude that
$f(B(x,\eps))\in\CU(F)$.

(a) $\Ra$ (b) and (c) $\Ra$ (d) are trivial.\finesp

If $A$ is a closed convex balanced, bounded subset of $F$, $F_A$ will
denote the Banach space obtained by taking the linear span of $A$
with the norm given by its Minkowski functional.

\begin{Th}
Let $\mathcal U$ be a closed surjective operator ideal, and
$f\in\HEF$. The following assertions are equivalent:

{\rm (a)} $f\in\HUEF$;

{\rm (b)} there is a closed convex, balanced subset $K\in \CU(F)$
such that $f$ is a holomorphic mapping from $E$ into $F_K$;

{\rm (c)} there is a Banach space $G$, a mapping $g\in{\mathcal
H}(E,G)$ and an operator $T\in{\mathcal U}(G,F)$ such that
$f=T\circ g$.
\end{Th}

\Proof
(a) $\Ra$ (b) follows the ideas in the proof of
\cite[Proposition~3.5]{AS} and \cite[Theorem~3.7]{RyWCo}.

For each $m\in\N$ and $x\in E$, define
$$
A_m(x):= \left\{ \lambda y:y\in B\left( x,\frac{1}{m}\right)
\mbox{ and } |\lambda |\leq 1\right\}
$$
and
$$
U_m:=\bigcup \left\{ B\left( x,\frac{1}{m}\right) : \|x\|\leq m
\mbox{ and } \|f\|_{A_m(x)}\leq m\right\} \, .
$$
For each $x\in E$ there is a neighbourhood of the compact set
$\{\lambda x:|\lambda |\leq 1\}$ on which $f$ is bounded.
Hence, there is $m\in\N$ so that $\|f\|_{A_m(x)}\leq m$, which shows
that $E=\cup_{m=1}^\infty U_m$.

Let $W_m$ be the balanced hull of $U_m$. Since the sets $A_m(x)$
are balanced, we have $|f(x)|\leq m$ for all $x\in W_m$.
Let $V_m:=2^{-1}W_m$.
We have $E=\cup_{m=1}^\infty V_m$ and hence
\begin{eqnarray}
\label{fE}
f(E)=\bigcup_{m=1}^\infty f(V_m)\, .
\end{eqnarray}
For each $k,m\in\N$, define
$$
K_{mk}:=\bar{\Gamma}\left( \frac{1}{k!}\, d^k\!f(0)(W_m)\right)
\in \CU(F)\, .
$$
By Proposition~\ref{propCU}, we obtain that the set
$$
K_m:=\left\{\sum_{k=0}^\infty 2^{-k}z_k:z_k\in K_{mk}\right\}
$$
belongs to $\CU(F)$. Easily, $f(V_m)\subseteq K_m$. Hence
$f(V_m)\in\CU(F)$ for all $m\in\N$. By Proposition~\ref{propCU},
we can select numbers $\beta_m>0$ with $\sum\beta_m<\infty$ so
that
$$
K:=\bar{\Gamma}\left( \bigcup_{m=1}^\infty \beta_mf(V_m)\right)
\in\CU(F)\, .
$$
It follows from (\ref{fE}) that $f$ maps $E$ into $F_K$.

It remains to show that $f\in{\mathcal H}(E,F_K)$. Let $x\in E$.
Easily, there are $\eps>0$ and $r\in\N$  such that $f(B(x,2\eps
))\subseteq rK$. By Lemma~\ref{Cauchy},
\begin{eqnarray}
\label{contmK}
\frac{1}{k!}\, d^k\!f(x)\left( B(0,2\eps )\right)\subseteq rK
\end{eqnarray}
for all $k\in\N$.
Now, for each $n\in \N$ and $a\in B(0,\eps )$, we have
$$
f(x+a)-\sum_{k=0}^n\frac{1}{k!}\, d^k\!f(x)(a) =
2^{-n}\sum_{k=n+1}^\infty 2^{n-k}\frac{1}{k!}\, d^k\!f(x)(2a)\, .
$$
Since $K$ is convex and closed, we get from (\ref{contmK}) that
$$
\sum_{k=n+1}^\infty 2^{n-k}\frac{1}{k!}\, d^k\!f(x)(2a)\in rK\, .
$$
Hence,
$$
f(x+a)-\sum_{k=0}^n\frac{1}{k!}\, d^k\!f(x)(a)\in 2^{-n}rK\, ,
$$
and so, the $F_K$-norm of the left hand side is less than or equal
to $2^{-n}r$, for all $a\in B(0,\eps )$. Thus, $f$ is holomorphic.

(b) $\Ra$ (c). It is enough to note that, by Lemma~\ref{surjid},
the natural inclusion $F_K\ra F$ belongs to ${\mathcal U}$.

(c) $\Ra$ (a). Each $x\in E$ has a neighbourhood $V_x$ such that
$g(V_x)$ is bounded in $G$. Hence, $f(V_x) =T(g(V_x))\in \CU (F)$.

\begin{Th}
\label{factor}
Let ${\mathcal U}$ be a closed surjective, factorizable operator
ideal and take a mapping $f\in\HEF$. Then $f \in\HUEF$ if and only
if there are a Banach space $G$, a mapping $g\in{\mathcal H}(E,G)$
and $T\in{\mathcal U}(G,F)$ such that $I_G\in{\mathcal U}$ and
$f=T\circ g$.
\end{Th}

\begin{Rem}{\rm
Theorem~\ref{factor} implies that, if ${\mathcal U}$ is the ideal
of weakly compact (resp.\ Rosenthal, Banach-Saks or Asplund)
operators and $f\in\HUEF$, then $f$ factors through a Banach space
$G$ which is reflexive (resp.\ contains a copy of $\ell_1$, has
the Banach-Saks property or is Asplund).

Moreover, if ${\mathcal U}=\{ T:T^q \mbox{ has separable
range}\}$, then $G$ is isomorphic to $G_1\times G_2$, with
$G_1^{**}$ separable and $G_2$ reflexive \cite{V}. If ${\mathcal
U}=\{ T:T^*\mbox{ is Rosenthal}\}$, then $G$ contains no copy of
$\ell_1$ and no quotient isomorphic to $c_0$ \cite{GO}.}
\end{Rem}

\end{document}